\documentclass[]{article}

\input epsf

\def\dsp{\displaystyle}
\def\pref#1{(\ref{#1})}

\def\bigO{{\cal O}}
\def\C{{\cal C}}

\def\FG#1#2#3#4{
{}_2F_1\left(
\begin{array}{c}
\begin{array}{cc} \hskip-10pt#1,{\ } #2 \end{array}\\
\begin{array}{c} \hskip-10pt#3 \end{array}
\end{array}
\hskip-8pt;\,#4
\right)}

\def\CHG#1#2#3{
{}_1F_1\left(
\begin{array}{c}
\begin{array}{c} \hskip-10pt#1 \end{array}\\
\begin{array}{c} \hskip-10pt#2 \end{array}
\end{array}
\hskip-8pt;\,#3
\right)}

\def\FFF#1#2#3#4#5#6{
{}_3F_2\left(
\begin{array}{c}
\begin{array}{c}\hskip-10pt#1,#2,#3\end{array}\\
\begin{array}{c}\hskip-10pt #4,#5\end{array}
\end{array}
\hskip-8pt;\,#6
\right)}

\def\sfrac#1#2{{{\lower.6ex
\hbox{$\scriptstyle#1$}}\over 
{\raise.7ex
\hbox{$\scriptstyle#2$}}}}

\begin{document}
 \title{
Asymptotics of a ${}_3F_2$ polynomial associated with the 
Catalan-Larcombe-French sequence}

   \author{Nico M. Temme\\       
     CWI,  P.O. Box 94079, NL-1090 GB Amsterdam,  The Netherlands\\           
   e-mail: {\tt  Nico.Temme@cwi.nl}   
   }

    \maketitle                                  
    \begin{abstract}                            
   \noindent                                    
   The large $n$ behaviour of the hypergeometric polynomial
   $$
   \FFF{-n}{\sfrac12}{\sfrac12}{\sfrac12-n}{\sfrac12-n}{-1}
   $$
   is considered by using integral representations of this polynomial.
   This ${}_3F_2$ polynomial is associated with the
   Catalan-Larcombe-French sequence. Several other representations 
   are mentioned, with references to the literature,
   and another asymptotic method is described by using 
   a generating function of the sequence. The results are similar to 
   those obtained by Clark (2004) who used a binomial sum for 
   obtaining an asymptotic expansion.

   \end{abstract}                               
                                                
   \vskip 0.8cm \noindent                       
   Mathematics Subject Classification 2000:     
   41A60, 33C20, 11B83, 33C10.                         
   \par\noindent                                
   Keywords \& Phrases:                         
   Catalan-Larcombe-French sequence, 
   asymptotic expansion, 
   hypergeometric polynomial,  
   modified Bessel function.

\section{The problem}
\label{sec.TP}
Find the large $n$ asymptotics of
\begin{equation}\label{P1}
f(n)=
\FFF{-n}{\sfrac12}{\sfrac12}{\sfrac12-n}{\sfrac12-n}{-1}
\end{equation}
Peter Larcombe conjectured  that
$\lim_{n\to\infty} f(n)=2$
and Tom Koornwinder gave a proof, based on dominated convergence. See 
for details
of the proof
\cite{Larcombe:2006:FPL}, where also a different representation of $f(n)$ is 
considered in the form 
\begin{equation}\label{P2}
f(n)=
2^{n}\FFF{-n}{-\sfrac12n}{\sfrac12-\sfrac12n}{\sfrac12-n}{\sfrac12-n}{1}.
\end{equation}
The equivalence of these two forms follows from a quadratic 
transformation of the ${}_{3}F_{2}-$functions as given in \cite[Ex.~4(iv), p.97]{Bailey:1964:GHS}, that is,
\renewcommand{\arraystretch}{1.5}
\begin{equation}\label{P3}\begin{array}{l}
\FFF{a}{b}{c}{1+a-b}{1+c-c}{z}=\\
\quad\quad\quad\quad
(1-z)^{-a}\FFF{\sfrac12a}{\sfrac12+\sfrac12a}{1+a-b-c}{1+a-b}{1+a-c}{\dsp{\frac{-4z}{(1-z)^{2}}}}.
\end{array}
\end{equation}
\renewcommand{\arraystretch}{1.0}
with $a=-n$, $b=c=\frac12$, and $z=-1$.
Another form is given by (see \cite[Eq. (A2)]{Larcombe:2004:ANG})
\begin{equation}
    \label{P4}
f(n)=\frac{n!}{2^{n}(\frac12)_{n}}\FFF{-n}{-n}{\frac12}{1}{\frac12-n}{-1}.
\end{equation}

In \cite{Clark:2004:AEC} an asymptotic expansion of $\frac12f(n)$ has been 
derived. The asymptotic analysis is based on the representation
\begin{equation}
    \label{P5}
P_{n}=\frac{1}{n!}\sum_{p+q=n}{2n\choose p}
{2q\choose q}\frac{(2p)!\,(2q)!}{p!\,q!}.
\end{equation}
By using the relation
\begin{equation}
    \label{P6}
(2n)!=2^{2n}n!(\sfrac12)_{n},\quad n=0,1,2,\ldots,\end{equation}
it is straightforward to verify that  \pref{P5} can be written 
as 
\begin{equation}
    \label{P7}
P_{n}=\frac{2^{4n}}{n!}\sum_{p=0}^{n}\frac{(\frac12)_{p}(\frac12)_{p}(\frac12)_{n-p}(\frac12)_{n-p}}{p!\,(n-p)!}.
\end{equation}
By using
\begin{equation}
    \label{P8}
(a)_{n-k}=(-1)^{k}\frac{(a)_{n}}{(1-a-n)_{k}},
\end{equation}
it follows that 
\begin{equation}
    \label{P9}
P_{n}=\frac{2^{4n}(\frac12)_{n}(\frac12)_{n}}{n!\,n!}
\sum_{p=0}^{n}(-1)^{p}\frac{(-n)_{p}(\frac12)_{p}(\frac12)_{p}}{p!\,(\frac12-n)_{p}(\frac12-n)_{p}},
\end{equation}
that is,
\begin{equation}
\label{P10}
P_{n}=\frac{2^{4n}(\frac12)_{n}(\frac12)_{n}}{n!\,n!}
\FFF{-n}{\sfrac12}{\sfrac12}{\sfrac12-n}{\sfrac12-n}{-1},
\end{equation}
which gives the relation with $f(n)$ by using \pref{P1}:
\begin{equation}
\label{P11}
P_{n}=\frac{2^{4n}(\frac12)_{n}(\frac12)_{n}}{n!\,n!}\,f(n)={2n 
\choose n}^{2}\,f(n).
\end{equation}

The numbers $P_{n}$ are for $n=0,1,2,,\ldots$ known as the elements of the 
sequence (A053175) $\{1,8,80, 896,10816,\ldots\}$, called the {\em 
Catalan-Larcombe-French} sequence, which is originally discussed by 
Catalan \cite{Catalan:1887:SNS}. See the 
{\em On-Line Encyclopedia of Integer Sequences}
http://www.research.att.com/~njas/sequences/ .

In this paper we derive a complete asymptotic expansion of the 
numbers $P_{n}$ by using integral representations of the 
corresponding ${}_{3}F_{2}-$functions. Our results are the same as 
those obtained by Clark \cite{Clark:2004:AEC}, who used the binomial 
sum in \pref{P5} without reference to the ${}_{3}F_{2}-$functions. 
\section{Transformations}
\label{sec.TR}
We derive an integral representation of the ${}_{3}F_{2}-$function of
\pref{P1} by using several transformations for special functions. We 
start with the beta integral
\begin{equation}\label{T0}
B(x,y)=\frac{\Gamma(x)\Gamma(y)}{\Gamma(x+y)}=\int_{0}^{1}t^{x-1}(1-t)^{y-1}\,dt
\end{equation}
and use it in the form
\begin{equation}\label{T1}
\frac{\left(\sfrac12\right)_k}{\left(\sfrac12-n\right)_k}=
\frac{(-1)^kn!}{\sqrt{{\pi}}\Gamma(n+\sfrac12)}\int_0^1
t^{k-\frac12}(1-t)^{n-k-\frac12}\,dt,\quad k=0,1,\ldots, n.
\end{equation}
We substitute this in the representation of the 
${}_{3}F_{2}-$function in \pref{P1}
\begin{equation}\label{T1a}
\FFF{-n}{\sfrac12}{\sfrac12}{\sfrac12-n}{\sfrac12-n}{-1}=
\sum_{k=0}^{n}(-1)^{k}\frac{(-n)_{k}(\frac12)_{k}(\frac12)_{k}}{k!\,(\frac12-n)_{k}(\frac12-n)_{k}}.
\end{equation}
This gives after performing the $k-$summation
\begin{equation}\label{T2}
f(n)=\frac{n!}{\pi\,(\sfrac12)_n}\int_0^1t^{-\frac12}(1-t)^{n-\frac12}
\FG{-n}{\sfrac12}{\sfrac12-n}{\frac{t}{1-t}}\,dt.
\end{equation}
We substitute $t=\sin^{2}(\theta/2)$ and obtain
\begin{equation}\label{T3}
f(n)=\frac{n!}{\pi\,(\sfrac12)_n}\int_0^{\pi}\cos^{2n}(\theta/2)
\FG{-n}{\sfrac12}{\sfrac12-n}{\tan^{2}(\theta/2)}\,d\theta.
\end{equation}
We apply a quadratic transformation (see \cite[Eq.~15.3.26]{Abramowitz:1964:HMF}) to obtain
\begin{equation}\label{T4}
f(n)=\frac{n!}{\pi\,(\sfrac12)_n}\int_0^{\pi}
\FG{-\sfrac12n}{\sfrac12-\sfrac12n}{\sfrac12-n}{\sin^{2}\theta}\,d\theta,
\end{equation}
and use the representation of the Legendre polynomial
\begin{equation}\label{T5}
P_{n}(x)=\frac{(2n)!}{2^{n}n!\,n!}x^{n}
\FG{-\sfrac12n}{\sfrac12-\sfrac12n}{\sfrac12-n}{x^{-2}}.
\end{equation}
This follows from \cite[Eq.~(22.3.8)]{Abramowitz:1964:HMF} and gives
\begin{equation}\label{T6}
f(n)=\frac{2^{-n}n!\,n!}{\pi\,(\sfrac12)_n\,(\sfrac12)_n}\int_0^\pi
\sin^n\theta\  P_n\left(\frac1{\sin\theta}\right)\,d\theta.
\end{equation}
Next, consider (see \cite[p.~204]{Temme:1996:SFI})
\begin{equation}\label{T7}
P_{n}(z)=\frac1{\pi}\int_{0}^{\pi}
\left(z+\sqrt{z^{2}-1}\cos\psi\right)^{n}\,d\psi, \quad 
n=0,1,2,\ldots,
\end{equation}
which gives the double integral
\begin{equation}\label{T8}
f(n)=\frac{n!\,n!}{\pi^{2}\,(\sfrac12)_n\,(\sfrac12)_n}\int_0^\pi
\int_0^\pi\left(\frac{1+\cos\theta\cos\psi}{2}\right)^n\,d\theta\,d\psi.
\end{equation}

   \begin{figure}                                                           
   \begin{center}                                                           
   \epsfxsize=6cm \epsfbox{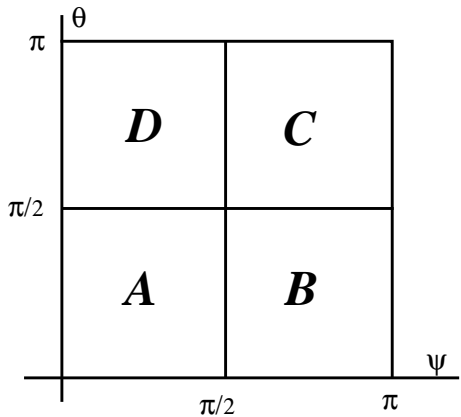}                                        
   \end{center}                                                             
   \caption{\small The domain of integration of the integral in \pref{T8}   
    and subdomains  $A$, $B$, $C$ and $D$. \label{larc1}}                    
   \end{figure}                                                             

\section{Asymptotic analysis}
\label{sec.AA}

The landscape of the integrand in \pref{T8} shows peaks at the boundary points
$(0,0)$ and $(\pi,\pi)$, where it assumes the value $1$. Along 
the interior lines $\theta=\frac12\pi$ and $\psi=\frac12\pi$ the 
integrand has the value $2^{-n}$. Inside the squares $A$ and $C$, see 
Figure~\ref{larc1}, the value of the integrand is between $2^{-n}$ and 
$1$, in the squares $B$ and $D$ it is between $0$ and 
$2^{-n}$. In addition, the contributions from $A$ and $C$ are the same, and 
also those from $B$ and $D$ are the same.

From an asymptotic point of view it follows that the integral over 
the full square equals twice the integral over $A$, with an error 
that is of order $\bigO(2^{-n})$, while the total integral is of order 
$\bigO(1)$, as $n$ is large. Hence, we concentrate on the integral 
over $A$, and write for large values of $n$
\begin{equation}\label{A1}
f(n)=2\frac{2^{-n}n!\,n!}{\pi^{2}\,(\sfrac12)_n\,(\sfrac12)_n}
\left[
\int_0^{\frac12\pi}
\int_0^{\frac12\pi}(1+\cos\theta\cos\psi)^n\,d\theta\,d\psi+E_{n}\right],
\end{equation}
where $E_{n}=\bigO(2^{-n})$. Next, we neglect $E_{n}$ and 
put $u=\sin(\theta/2), v=\sin(\psi/2)$, and obtain
\begin{equation}\label{A2}
f(n)\sim\frac{8\,n!\,n!}{\pi^{2}\,(\sfrac12)_n\,(\sfrac12)_n}
\int_0^{\frac12\sqrt{2}}
\int_0^{\frac12\sqrt{2}}
\left(1-u^{2}-v^{2}+
2u^{2}v^{2}\right)^n\,\frac{du}{\sqrt{1-u^{2}}}\,\frac{dv}{\sqrt{1-v^{2}}}.
\end{equation}

   For the integrals in \pref{A1} and \pref{A2} asymptotic expansions
can be obtained by using Laplace's method for double integrals; see 
\cite[\S~VIII.10]{Wong:2001:AAI}). In our case a simpler 
approach is based on neglecting a part of square $A$ by introducing polar coordinates
\begin{equation}\label{A3}
u=r\cos t, \quad v=r\sin t,\quad 0\le r\le \sfrac12\sqrt{2}, \quad 0\le t\le 
\sfrac12\pi.
\end{equation}
This gives (again we make an error in the integral that is of order $\bigO(2^{-n})$)
\begin{equation}\label{A4}
f(n)\sim\frac{8\,n!\,n!}{\pi^{2}\,(\sfrac12)_n\,(\sfrac12)_n}\int_0^{\pi/2}
\int_0^{\frac12\sqrt{2}}\frac{(1-r^{2}+
2r^{4}\cos^{2}t\sin^{2}t)^n r\,dr\,dt}
{\sqrt{(1-r^{2}\cos^{2}t)(1-r^{2}\sin^{2}t)}}.
\end{equation}
We change $r^{2}$ into $r$, and obtain
\begin{equation}\label{A5}
f(n)\sim\frac{4\,n!\,n!}{\pi^{2}\,(\sfrac12)_n\,(\sfrac12)_n}\int_0^{\pi/2}
\int_0^{\frac12}\frac{(1-r+
2r^{2}\cos^{2}t\sin^{2}t)^n \,dr\,dt}
{\sqrt{(1-r\cos^{2}t)(1-r\sin^{2}t)}}.
\end{equation}

First the standard method for obtaining asymptotic expansions of a 
Laplace-type integral can be used (for the $r-$integral). The second 
step is done by integrating the coefficients of this expansion with 
respect to $t$.

For the $r-$integral we transform the variable of integration by putting
\begin{equation}\label{A6}
w=-\ln\left(1-r+2r^{2}\cos^{2}t\sin^{2}t\right).
\end{equation}
This mapping is one-to-one for  
$r\in[0,\frac12]$, uniformly with respect to $t\in[0,\frac12\pi]$,
with corresponding $w-$interval $[0,w_{0}]$, where 
$w_{0}=w(\frac12)$.

We obtain
\begin{equation}\label{A7}
f(n)\sim\frac{4\,n!\,n!}{\pi^{2}\,(\sfrac12)_n\,(\sfrac12)_n}\int_0^{\pi/2}
\int_0^{w_{0}}e^{-nw}F(w,t)\,dw\,dt,
\end{equation}
where
\begin{equation}\label{A8}
F(w,t)=\frac{1}{\sqrt{(1-r\cos^{2}t)(1-r\sin^{2}t)}}\,\frac{dr}{dw}.
\end{equation}
\section{Asymptotic expansion}
\label{sec.AE}
We obtain the asymptotic expansion of $w-$integral in 
\pref{A7} by using Watson's lemma (see \cite[\S~I.5]{Wong:2001:AAI}).

The function $F(w,t)$ is  analytic in a 
neighborhood of the origin of the $w-$plane. We expand
\begin{equation}\label{E1}
F(w,t)=\sum_{k=0}^{\infty} c_{k}(t)w^{k}
\end{equation}
and substitute this expansion in \pref{A7}. Interchanging the order 
of summation and integration, and replacing the interval of the $w-$integrals
by $[0,\infty)$ (a standard procedure in asymptotics) we obtain
\begin{equation}\label{E2}
f(n)\sim\frac{4\,n!\,n!}{n\pi^{2}\,(\sfrac12)_n\,(\sfrac12)_n}
\sum_{k=0}^{\infty}C_{k}\frac{k!}{n^{k}},\quad n\to\infty,
\end{equation}
where 
\begin{equation}\label{E3}
C_{k}=\int_{0}^{\pi/2}c_{k}(t)\,dt,\quad k=0,1,2,\ldots\,.
\end{equation}

The coefficients $c_{k}(t)$ can be obtained by the following method.
First we need the inverse of the transformation defined in \pref{A6}. 
That is, we need coefficients $b_{k}$ in the expansion
\begin{equation}\label{E3a}
r(w)=\sum_{k=0}^{\infty}b_{k}(t)w^{k}.
\end{equation}
We can find $r(w)$ from \pref{A6} as a solution of a quadratic 
equation, with the condition $r(w)\sim w$ as $w\to 0$, that is, 
$b_{0}(t)=1$. However, we 
can also differentiate \pref{A6} with respect to $r$ and substitute 
the expansion \pref{E3a}, and solve for the coefficients 
$b_{k}(t)$. When we have these coefficients we can expand $F(w,t)$ of
\pref{A8} and find $c_{k}(t)$.

The first few coefficients $c_{k}(t)$ are
\renewcommand{\arraystretch}{1.5}
\begin{equation}
\label{E4}
\begin{array}{lll}
c_{0}(t) &= & 1,\\
c_{1}(t) &= & \frac12( -1+8s^2-8s^4),\\
c_{2}(t) &= & \frac18(1-28s^2+220s^4-384s^6+192s^8),\\
c_{3}(t) &= & \frac1{48}(-1+92s^2-1628s^4+10752s^6-24576s^8+23040s^{10}-7680s^{12}),\\
c_{4}(t) &= & 
\frac1{384}(1-280s^2+10024s^4-130848s^6+773904s^8-2054400s^{10}+\\
& & 2691840s^{12}-1720320s^{14}+430080s^{16}),\\
c_{5}(t) &= & 
\frac1{3840}(-1+848s^2-55328s^4+1259040s^6-13396560s^8+\\
& & 73983360s^{10}-215329920s^{12}+349224960s^{14}-\\
& & 319549440s^{16}+154828800s^{18}-30965760s^{20})
\end{array}
\end{equation}
\renewcommand{\arraystretch}{1.0}
where $s=\sin^{2}t$. For the corresponding $C_{k}$ we have
\begin{equation}
\label{E5}
C_{0}=\sfrac{1}{2}\pi,\quad
C_{1}=0,\quad
C_{2}=\sfrac{1}{8}\pi,\quad
C_{3}=\sfrac{1}{8}\pi,\quad
C_{4}=\sfrac{55}{384}\pi,\quad
C_{5}=\sfrac{11}{64}\pi.
\end{equation}

As a next step we can replace in \pref{E2} the ratios 
$n!/(\frac12)_n$ by the asymptotic expansion
\begin{equation}
\label{E6}
\frac{n!}{(\frac12)_n}=\sqrt{\pi}\frac{\Gamma(n+1)}{\Gamma(n+\frac12)}\sim 
\sqrt{\pi n}\sum_{k=0}^{\infty}\frac{\gamma_{k}}{n^{k}},
\end{equation}
where
\begin{equation}
\label{E7}
\gamma_{0}=1,\quad
\gamma_{1}=\sfrac{1}{8},\quad
\gamma_{2}=\sfrac{1}{128},\quad
\gamma_{3}=-\sfrac{5}{1024},\quad
\gamma_{4}=-\sfrac{21}{32768},\quad
\gamma_{5}=\sfrac{399}{262144}.
\end{equation}

This finally gives
\begin{equation}
\label{E8}
f(n)\sim 
2\left(1+\frac1{4n}+\frac{17}{32n^{2}}+\frac{207}{128n^{3}}+\frac{14875}{2048n^{4}}+
\frac{352375}{8192n^{5}}+\ldots\right).
\end{equation}

\section{An alternative method}
\label{sec.AM}
The numbers $P_{n}$ were proposed as ``Catalan'' numbers by an associate of
Catalan.  They appear as coefficients in the series expansion of an
elliptic integral of the first kind
\begin{equation}
\label{M1}
K(k)= \int_{0}^{\frac12\pi}\frac{1} {\sqrt{1 -k^2\sin^2t}}\,dt,
\end{equation}
which is transformed and
written as a power series in $k$ (through an intermediate variable); 
this gives a generating function for the sequence $\{P_{n}\}$.
For details we refer to \cite{Jarvis:2004:LRB}.

In \cite{Temme:2003:LPC} a generating function for the numbers $P_{n}$
is given in terms of the square of a
modified Bessel function, and we use this approach to obtain an asymptotic 
expansion of $f(n)$.
See also \cite{Larcombe:2004:ANG} for details on this generating 
function.

We consider numbers $F_{n}$ defined as coefficients in the generating 
function
\begin{equation}
    \label{M2}
\left[e^{w/2}\,I_0(w/2)\right]^2 = \sum_{n=0}^\infty F_n w^n.
\end{equation}
By considering the relation of the Bessel function with the confluent 
hypergeometric functions (see \cite[Eq.~13.6.3]{Abramowitz:1964:HMF}),
\begin{equation}
    \label{M3}
e^{z}I_{\nu}(z)=\frac{(\frac12z)^{\nu}}{\Gamma(\nu+1)}\CHG{\nu+\frac12}{2\nu+1}{2z},
\end{equation}
we can write \pref{M2} in the form (see also \cite[Eq.~13.1.27]{Abramowitz:1964:HMF}),
\begin{equation}
    \label{M4}
\left[\CHG{\frac12}{1}{w}\right]^{2}=e^{2w}\left[\CHG{\frac12}{1}{-w}\right]^{2} = \sum_{n=0}^\infty F_n w^n.
\end{equation}
This gives the representation for $F_{n}$:
\begin{equation}
    \label{M5}
F_{n}=\sum_{k=0}^{n}\frac{(\frac12)_{k}}{k!\,k!}\,\frac{(\frac12)_{n-k}}{(n-k)!\,(n-k)!}.
\end{equation}
By using \pref{P8} it follows that
\begin{equation}
    \label{M7}
F_{n}=\frac{(\frac12)_{n}}{n!\,n!}\sum_{k=0}^{n}(-1)^{k}\frac{(-n)_{k}(-n)_{k}(\frac12)_{k}}{(\frac12-n)_{k}k!\,k!},
\end{equation}
or
\begin{equation}
    \label{M8}
F_{n}=\frac{(\frac12)_{n}}{n!\,n!}\FFF{-n}{-n}{\frac12}{1}{\frac12-n}{-1}.
\end{equation}
It follows from \pref{P4} that
\begin{equation}
    \label{M9}
f(n)=\frac{n!\,n!\,n!}{2^n\,(\sfrac12)_n\,(\sfrac12)_n} \ F_n.
\end{equation}

From \pref{M2} we obtain
\begin{equation}  
    \label{M10}
F_{n}=\frac1{2\pi i}\int_{\C}\frac{\left[e^{w/2}\,I_0(w/2)\right]^2}{w^{n+1}}\, dw
=\frac1{2\pi i}\int_{\C}\frac{e^{2w}}{w^{n+1}}\,h(w)\, dw,
\end{equation}
where
\begin{equation}
    \label{M11}
h(w)=\left[e^{-w/2}\,I_0(w/2)\right]^2=\left[\CHG{\frac12}{1}{-w}\right]^{2},
\end{equation}
and the contour $\C$ is a circle around the
origin, or any contour that can be obtained from this circle by using 
Cauchy's theorem. The main contribution comes from the saddle point of 
$\frac{e^{2w}}{w^{n+1}}$, that is from $w=w_0=n/2$.

In the standard saddle point method (see 
\cite[\S~II.4]{Wong:2001:AAI}) a quadratic transformation is used to 
bring the main part of the integrand in the form of a Gaussian. We 
can obtain the same expansion by just expanding the function $h(w)$ 
(which is slowly varying for $w>0$)
at the saddle point. 

First we expand (see \cite[Eq.~13.4.9]{Abramowitz:1964:HMF})
\begin{equation}
    \label{M12}
\CHG{\frac12}{1}{-w}=\sum_{k=0}^{\infty}a_{k}(w-w_{0})^{k},\quad 
a_{k}=\frac{(-1)^{k}(\frac12)_{k}}{k!\,k!}\CHG{\frac12+k}{1+k}{-w_{0}}
\end{equation}
and next
\begin{equation}
    \label{M13}
h(w)=\sum_{k=0}^{\infty}A_{k}(w-w_{0})^{k}.
\end{equation}
We substitute this expansion in the second integral in \pref{M10} 
and obtain the convergent expansion
\begin{equation}
    \label{M14}
F_{n}=\sum_{k=0}^{\infty}A_{k}\Phi_{k},\quad \Phi_{k}=\frac1{2\pi 
i}\int_{\C}\frac{e^{2w}}{w^{n+1}}\,(w-w_{0})^{k}\, dw.
\end{equation}
The functions $\Phi_{k}$ can be evaluated by using the recursion 
formula
(which easily follows from integrating by parts)
\begin{equation}
    \label{M15}
\Phi_{k}=-\sfrac12(k-1)(\Phi_{k-1}+w_{0}\Phi_{k-2}),\quad 
\Phi_{0}=\frac{2^{n}}{n!},\quad \Phi_{1}=0.
\end{equation}

An asymptotic expansion can be obtained by using a well-known 
expansion for
$a_{k}$ defined in \pref{M12}. We have (as follows from 
\cite[Eq.~13.5.1]{Abramowitz:1964:HMF})
\begin{equation}
    \label{M16}
\CHG{a}{c}{-x}\sim 
x^{-a}\frac{\Gamma(c)}{\Gamma(c-a)}\sum_{m=0}^{\infty}
\frac{(a)_{m}(1+a-c)_{m}}{m!\,x^{m}},\quad x\to+\infty,
\end{equation}
from which we can obtain expansions for $a_{k}$ and $A_{k}$ for large 
values of $w_{0}=n/2$. By using these expansions in \pref{M14} we 
obtain an expansion for $F_{n}$, and finally for $f(n$ by using 
\pref{M9}. This expansion is the same as the one in \pref{E8}.

\section*{Acknowledgment}
I wish to thank Peter Larcombe for suggesting this problem, 
for encouraging me to investigate the asymptotic properties of the 
sequence $\{P_{n}\}$, and for introducing me to the literature, in 
particular to the papers [4] -- [7].

\bibliographystyle{plain}


\end{document}